\documentclass [12pt] {amsart}
\usepackage{amsthm,amsmath,amssymb,amscd,graphics,enumerate}
\usepackage{latexsym}
\usepackage{verbatim}
\usepackage{enumerate}
\usepackage{amsthm}
\usepackage{amsmath}
\usepackage{amscd}

\newcounter{theorem}[section]
\numberwithin{equation}{section}

\newtheorem{Thm}[theorem]{Theorem}
\newtheorem{Conj}[theorem]{Conjecture}
\newtheorem{Situ}[theorem]{Situation}
\newtheorem{cor}[theorem]{Corollary}

\newtheorem{lemma}[subsection]{Lemma}

\theoremstyle{definition}

\theoremstyle{remark}

\theoremstyle{remark}

\numberwithin{equation}{section}
{\theoremstyle{remark}
 
 }

{\theoremstyle{definition}
 }

\newcommand{\Aut}{{\operatorname{Aut}}}

\newcommand{\mD}{\mathcal{D}}

\newcommand{\X}{\mathcal{X}}
\newcommand{\Y}{\mathcal{Y}}

\newcommand{\D}{\mathcal{D}}

\newcommand{\on}{\operatorname}
\newcommand{\mP}{\mathcal{P}}

\def\<{\left\langle}
\def\>{\right\rangle}

\begin{document}

\title[Orbifold GW in codimension one]{On orbifold Gromov-Witten theory in codimension one}

\author{Hsian-Hua Tseng}
\address{Department of Mathematics\\ Ohio State University\\ 100 Math Tower, 231 West 18th Ave.\\Columbus\\ OH 43210\\ USA}
\email{hhtseng@math.ohio-state.edu}

\author{Fenglong You}
\address{Department of Mathematics\\ Ohio State University\\ 100 Math Tower, 231 West 18th Ave.\\Columbus\\ OH 43210\\ USA}
\email{you.111@osu.edu}
\date{\today}

\begin{abstract}
We propose a conjectural determination of the Gromov-Witten theory of a root stack along a smooth divisor. We verify our conjecture under an additional assumption.
\end{abstract}
\maketitle
\section{Introduction}
We work over $\mathbb{C}$. 

A folklore result, proven only recently in \cite{gs}, states that for a smooth Deligne-Mumford stack, the stack structures in codimension one are obtained by applications of root constructions along divisors, as defined in \cite{c} and \cite{agv}. From the perspective of Gromov-Witten theory, this motivates the attempts to understand Gromov-Witten theory of root stacks along divisors.

Let $X$ be a smooth proper Deligne-Mumford stack with projective coarse moduli space. Let $D\subset X$ be an irreducible smooth divisor. For $r\in \mathbb{N}$, let $\X_r$ denote the stack of $r$-th roots along $D$. There is a natural map $\X_r\to X$. Let $\D_r\subset \X_r$ be the divisor lying over $D\subset X$. 

We propose the following

\begin{Conj}\label{conj:gw_root}
The Gromov-Witten theory of $\X_r$ is determined by the Gromov-Witten theory of $X, D$, and the restriction map $H^*(IX)\to H^*(ID)$. 
\end{Conj}

Consider the following 

\begin{Situ}\label{nonstacky}
Let $X$, $D$, and $r$ be as above. In addition, assume that $D$ is disjoint from the locus of stack structures of $X$.
\end{Situ}

The main result of this paper is the following evidence supporting Conjecture \ref{conj:gw_root}:

\begin{Thm}\label{thm:main}
Conjecture \ref{conj:gw_root} holds true in Situation \ref{nonstacky}.
\end{Thm}
A proof of Theorem \ref{thm:main} is given in Section \ref{sec:pf_thm}. A discussion on the more general case is given in Section \ref{subsec:gen}.

H. -H. T. is supported in part by a Simons Foundation Collaboration Grant.

\section{On Conjecture \ref{conj:gw_root}}
Our study of Gromov-Witten theory of $\X_r$ relies on relative Gromov-Witten theory of Deligne-Mumford stack pairs, as defined and studied in \cite{af}. We also make heavy use of the results in \cite{mp}.

\subsection{Proof of Theorem \ref{thm:main}}\label{sec:pf_thm}

\subsubsection{Degeneration}
The deformation to the normal cone of $\D_r\subset \X_r$ yields a degeneration of $\X_r$ to a nodal stack $\X_r\cup_{\D_0} \Y$, where $\Y:=\mathbb{P}(N\oplus \mathcal{O})$ is obtained from the normal bundle $N$ of $\D_r\subset \X_r$. The gluing is done by identifying the zero section $\D_0\subset \Y$ of $\Y\to \D_r$ with $\D_r\subset \X_r$. 

The degeneration formula \cite[Theorem 0.4.1]{af} applied\footnote{We should also invoke deformation invariance \cite[Theorem 4.6.1]{af}.} to this degeneration shows that Gromov-Witten invariants of $\X_r$ are expressed in terms of relative Gromov-Witten invariants of $(\X_r, \D_r)$ and of $(\Y, \D_0)$. Gromov-Witten theory of $(\X_r, \D_r)$ is solved in Section \ref{subsubsec:relGW}, and Gromov-Witten theory of $(\Y, \D_0)$ is solved in Section \ref{subsubsec:L_H}. The proof of Theorem \ref{thm:main} is complete.

\subsubsection{Relative theory}\label{subsubsec:relGW}
In Situation \ref{nonstacky}, the inertia stack $I\X_r$ is a disjoint union of the inertia stack $IX$ and $r-1$ $\mu_r$-gerbes over $D$. Marked points receiving insertions from these $\mu_r$ gerbes are constrained to map to $\D_r$, and hence are relative marked points. It follows that insertions for a Gromov-Witten invariant of $(\X_r, \D_r)$ can be expressed as pull-backs of classes from $IX$. By \cite[Proposition 4.5.1]{af}, Gromov-Witten invariants of $(\X_r, \D_r)$ are equal to Gromov-Witten invariants of $(X, D)$. 

By \cite[Theorem 2]{mp}, if $X$ is a scheme, Gromov-Witten invariants of $(X, D)$ are determined by Gromov-Witten invariants of $X, D$, and the restriction map $H^*(X)\to H^*(D)$. In Situation \ref{nonstacky}, $X$ need not be a scheme. However $D$ is a scheme, and the proof of \cite[Theorem 2]{mp} still applies. More precisely, the degeneration formula applied to the deformation to the normal cone of $D\subset X$ shows that Gromov-Witten invariants of $(X, D)$ are obtained as solutions of a system of linear equations whose coefficients are given by Gromov-Witten invariants of $X$ and $(\mathbb{P}(N_{D/X}\oplus\mathcal{O}), D_0)$. By   \cite[Theorem 1]{mp}, Gromov-Witten invariants of $(\mathbb{P}(N_{D/X}\oplus\mathcal{O}), D_0)$ are expressed in terms of Gromov-Witten invariants of $D$, and by \cite[Lemma 4]{mp}\footnote{The proof of this Lemma in our setting is identical to the one in \cite{mp}.} this system is nonsingular. Gromov-Witten invariants of $(X, D)$ are thus solved.

\subsubsection{Leray-Hirsch}\label{subsubsec:L_H}
The projection map $\Y\to \D_r$ gives $\Y$ the structure of a fiber bundle with fibers the toric stack $\mathbb{P}^1_{r,r}$. By construction $\Y$ admits a fiberwise $\mathbb{C}^*$-action. The $\mathbb{C}^*$-fixed loci are the zero section $\D_0$ and the infinity section $\D_\infty$, both of which are isomorphic to $\D_r$. Virtual localization expresses Gromov-Witten invariants of $\Y$ in terms of certain Hurwitz-Hodge integrals over stable map moduli for $\D_r$. The orbifold quantum Riemann-Roch formula \cite{Tseng} applies to express these Hurwitz-Hodge integrals in terms of Gromov-Witten invariants of $\D_r$. By the main results of \cite{AJT, AJT2}, Gromov-Witten invariants of the $\mu_r$-gerbe $\D_r$ over $D$ are expressed in terms of Gromov-Witten invariants of $D$. 

Applying relative virtual localization, we see that Gromov-Witten invariants of $(\Y, \D_0)$ are expressed in terms of integrals on stable map moduli to $\Y$ (which are solved in the aforementioned manner) and rubber invariants associated to $(\Y, \D_0\cup \D_\infty)$. These rubber invariants contain {\em target descendants}. In the same way as \cite[Section 1.5]{mp}, tautological equations associated to the target descendant can be applied to express these rubber invariants in terms of rubber invariants associated to $(\Y, \D_0\cup \D_\infty)$ {\em without} target descendants. The rigidification result \cite[Lemma 2]{mp}\footnote{The proof of this Lemma in our setting is identical to the one in \cite{mp}.} applies to express rubber invariants associated to $(\Y, \D_0\cup \D_\infty)$ without target descendants as Gromov-Witten invariants of $(\Y, \D_0\cup \D_\infty)$.

Therefore it remains to determine Gromov-Witten invariants of $(\Y, \D_0\cup \D_\infty)$. By construction, $\Y$ is obtained from the $\mathbb{P}^1$-bundle $Y:=\mathbb{P}(N_{D/X}\oplus\mathcal{O})$ over $D$ by applying the $r$-th root construction to the zero and infinity sections $D_0, D_\infty \subset Y$. It follows from \cite[Proposition 4.5.1]{af} that Gromov-Witten invariants of $(\Y, \D_0\cup \D_\infty)$ are equal to Gromov-Witten invariants of $(Y, D_0\cup D_\infty)$. By \cite[Theorem 1]{mp}, Gromov-Witten invariants of $(Y, D_0\cup D_\infty)$ are solved in terms of Gromov-Witten invariants of $D$. 

\subsubsection{Remarks}
In Situation \ref{nonstacky} the divisor $D$ is assumed to be a scheme. This assumption is restrictive: from the perspective of introducing stack structures in codimension one by root constructions, it is more natural to allow $D$ to have nontrivial stack structures. Nevertheless Situation \ref{nonstacky} does apply when one constructs {\em twisted curves} (in the sense of \cite{av}) by applying root constructions to nonsingular curves. Gromov-Witten theory of nonsingular curves is completely solved by \cite{op1, op2, op3}. Therefore we obtain the following:

\begin{cor}
Gromov-Witten theory of a nonsingular twisted curve is completely determined.
\end{cor}

It is reasonable to expect that the proof of Theorem \ref{thm:main} can be made into algorithms that compute Gromov-Witten invariants of $\X_r$. More quantitative results on Gromov-Witten theory of $\X_r$, such as close formulas for Gromov-Witten invariants in low genera, are also desirable. Results along this line will be pursued elsewhere.

\subsection{On the general case}\label{subsec:gen}
In order to address Conjecture \ref{conj:gw_root} beyond Situation \ref{nonstacky}, it is natural to try to use the strategy for the proof of Theorem \ref{thm:main}. To make this work, we need the following

\begin{Conj}\label{conj:stackyLH}
Let $D$ be a smooth proper Deligne-Mumford stack with projective coarse moduli space. Let $L$ be a line bundle on $D$. Put $Y:=\mathbb{P}(L\oplus \mathcal{O})$. Write $\pi: Y\to D$ for the natural projection and  $D_0, D_\infty\subset Y$ for the zero and infinity sections of $\pi$. 

Then the Gromov-Witten theory of $Y$, $(Y, D_0)$, and $(Y, D_0\cup D_\infty)$ are determined by the Gromov-Witten theory of $D$ and $c_1(L)\in H^2(D)$. 
\end{Conj}

Conjecture \ref{conj:stackyLH} is a generalization of \cite[Theorem 1]{mp} to orbifold Gromov-Witten theory. Arguments in Section \ref{subsubsec:L_H} can be understood as proving Conjecture \ref{conj:stackyLH} when $D$ is a root gerbe\footnote{See e.g. \cite{AJT, AJT2} for discussions on root gerbes.} over a scheme. It is also easy to see that the part of Conjecture \ref{conj:stackyLH} about Gromov-Witten theory of $Y$ holds true for any $D$: Simply apply virtual localization with respect to the fiber wise $\mathbb{C}^*$-action and orbifold quantum Riemann-Roch. 

The {\em relative} Gromov-Witten theories of $(Y, D_0)$ and $(Y, D_0\cup D_\infty)$ are much harder to solve. The difficulty already appears in solving Gromov-Witten invariants in {\em fiber classes}, i.e. curve classes $\beta\in H_2(Y)$ such that $\pi_*\beta=0$. A simple argument along the line of \cite[Section 1.2]{mp} shows that Gromov-Witten invariants in fiber classes are determined by relative Gromov-Witten invariants of fibers of $\pi:Y\to D$. For a point $x\in D$ with stabilizer group $G$, the fiber of $\pi$ over $x$ is of the form $[\mathbb{P}^1/G]$, where $G$ acts on $\mathbb{P}^1$ via its action on $\mathbb{C}^2=L_x\oplus \mathbb{C}$: $G$ acts on the fiber $L_x$ and acts trivially on $\mathbb{C}$. Gromov-Witten invariants of $[\mathbb{P}^1/G]$ {\em relative} to $[0/G]$ and/or $[\infty/G]$ are not known, and appear difficult to determine.

\appendix
\section{On relative orbifold Gromov-Wiiten invariants}
We discuss here some results on relative Gromov-Witten theory of Deligne-Mumford stacks. 

\subsection{Relative virtual localization}
Virtual localization for relative Gromov-Witten theory of scheme pairs is detailed in \cite{gv}. For Deligne-Mumford stack pairs, the assumption on global nonsingular embeddings needed in the proof of virtual localization in \cite{gp} may be verified by combining the constructions of nonsingular embeddings in \cite{gv} and \cite{agot}. Alternatively, one can invoke the proof of virtual localization in \cite{ckl} without nonsingular embedding.  

The virtual localization formula for relative Gromov-Witten invariants of Deligne-Mumford stack pairs may be derived in the same way as in \cite{gv}. The formula takes the same form as the one for scheme pairs. Stack structures at nodes result in additional scalar factors that need be carefully sorted out. However, the precise nature of these factors do not affect the qualitative application of relative virtual localization formula needed for the proof of Theorem \ref{thm:main}.

\subsection{Rubber invariants}
Let $D, Y, D_0, D_{\infty}$ be as in Conjecture \ref{conj:stackyLH}. Rubber invariants associated to $(Y, D_0\cup D_\infty)$, which are not defined in \cite{af}, may be defined in the same way as in e.g. \cite[Section 1.5]{mp}. 

In this Deligne-Mumford setting, rubber invariants with target descendants also satisfy divisor and dilaton equations and tautologial recursion relations. These equations take the same form\footnote{Note a typo in the dilaton equation in \cite[Section 1.5.4]{mp}: $\tau_1(1)$ should not appear on the right-hand side.} as their scheme-theoretic counterparts, see \cite[Section 1.5.4, Section 1.5.5]{mp}.

Rubber invariants without target descendants can be expressed in terms of Gromov-Witten invariants of $(Y, D_0\cup D_\infty)$. More precisely, \cite[Lemma 2]{mp} holds in this Deligne-Mumford setting, with the same proof.


\begin{thebibliography}{HLOY02}

\bibitem{af} D. Abramovich, B. Fantechi, {\em Orbifold techniques in degeneration formulas}, to appear in Annali della SNS, arXiv:1103.5132.


\bibitem{agot} D. Abramovich, T. Graber, M. Olsson, H.-H. Tseng, {\em On the global quotient structure of the space of twisted stable maps to a quotient stack}, J. Algebraic Geom. 16 (2007), no. 4, 731--751.

\bibitem{agv} D. Abramovich, T. Graber, A. Vistoli, {\em Gromov-Witten theory of Deligne-Mumford stacks}, Amer. J. Math. 130 (2008), no. 5, 1337--1398.


\bibitem{av} D Abramovich, A. Vistoli, {\em Compactifying the space of stable maps},  J. Amer. Math. Soc.  15  (2002),  no. 1, 27--75

\bibitem{AJT} E. Andreini, Y. Jiang, H.-H. Tseng,
{\em Gromov-Witten theory of root gerbes I: structure of genus $0$ moduli spaces,} J. Differential Geom. 99 (2015), no. 1, 1--45.

\bibitem{AJT2} E. Andreini, Y. Jiang, H.-H. Tseng, {\em Gromov-Witten theory of banded gerbes over schemes}, arXiv:1101.5996.

\bibitem{c} C. Cadman, {\em Using stacks to impose tangency conditions on curves}, Amer. J. Math. 129 (2007), no. 2, 405--427.

\bibitem{ckl} H.-L. Chang, Y.-H. Kiem, J. Li, {\em Torus localization and wall crossing for cosection localized virtual cycles}, arXiv:1502.00078.

\bibitem{gs} A. Geraschenko, M. Satriano, {\em A "bottom up" characterization of smooth Deligne-Mumford stacks}, arXiv:1503.05478.

\bibitem{gp} T. Graber, R. Pandharipande, {\em Localization of virtual classes}, Invent. Math. 135 (1999), no. 2, 487--518.
 
\bibitem{gv} T. Graber, R. Vakil, {\em Relative virtual localization and vanishing of tautological classes on moduli spaces of curves}, Duke Math. J. 130 (2005), no. 1, 1--37. 

\bibitem{mp} D. Maulik, R. Pandharipande, {\em A topological view of Gromov-Witten theory}, Topology 45 (2006), no. 5, 887--918.

\bibitem{op1} A. Okounkov, R. Pandharipande, {\em Gromov-Witten theory, Hurwitz theory, and completed cycles} Ann. of Math. (2) 163 (2006), no. 2, 517--560.

\bibitem{op2} A. Okounkov, R. Pandharipande, {\em The equivariant Gromov-Witten theory of $\mathbf{P}^1$}, Ann. of Math. (2) 163 (2006), no. 2, 561--605.

\bibitem{op3} A. Okounkov, R. Pandharipande, {\em Virasoro constraints for target curves}, Invent. Math. 163 (2006), no. 1, 47--108. 

\bibitem{Tseng} H.-H. Tseng,  {\em Orbifold quantum Riemann-Roch, Lefschetz, and Serre}, Geom. Topol. 14 (2010) 1--81.

\end{thebibliography}
\end{document}